\documentclass[12pt,draft]{amsart}
\usepackage{amsmath,amsthm,latexsym,amscd,amsbsy,amssymb}
\newcommand{\spep}{${\mathcal S}_{X} =
\{ X_{\alpha}, p_{\alpha}^{\beta}, A \}\;$}

\newcommand{\speq}{${\mathcal S}_{Y} =
\{ Y_{\alpha}, q_{\alpha}^{\beta}, A \}\;$}

\chardef\bslash=`\\ % p. 424, TeXbook

\makeatletter
\def\verbatim{\interlinepenalty\@M \@verbatim
  \leftskip\@totalleftmargin\advance\leftskip2pc
  \frenchspacing\@vobeyspaces \@xverbatim}
\makeatother
\newtheorem{thm}{Theorem}[section]
\newtheorem{cor}[thm]{Corollary}
\newtheorem{lem}[thm]{Lemma}
\newtheorem{pro}[thm]{Proposition}

\newtheorem*{A}{{\bf Characterization of injectives}}

\newtheorem*{C}{{\bf Retracts of uncountable
products of Banach spaces}}
\newtheorem*{D}{{\bf Retracts of uncountable
coproducts of Banach spaces}}

\theoremstyle{definition}
\newtheorem{defin}{Definition}[section]

\theoremstyle{remark}
\newtheorem{rem}{Remark}[section]

\numberwithin{equation}{section}

\begin{document}

%%%%%%% Begin Topmatter %%%%%%%%%%

\title[Continuous homomorphisms of Arens-Michael algebras]
{Continuous homomorphisms of Arens-Michael algebras}
\author{Alex Chigogidze}
\address{Department of Mathematics and Statistics,
University of Saskatche\-wan,
McLean Hall, 106 Wiggins Road, Saskatoon, SK, S7N 5E6,
Canada}
\email{chigogid@math.usask.ca}
\thanks{Author was partially supported by NSERC research grant.}


%%%%%%% End topmatter %%%%%%%%%

\begin{abstract}{It is shown (Theorem \ref{T:spectral})
that every continuous homomorphism of Arens-Michael algebras
can be obtained as the limit
of a morphism of certain projective systems consisting of 
Fr\'{e}chet algebras. Based on this we prove (Theorem
\ref{T:complementedpro}) that a complemented
subalgebra of an uncountable product of Fr\'{e}chet algebras is
topologically isomorphic to
the product of Fr\'{e}chet algebras.
These results are used to characterize (Theorem \ref{T:inj})
injective objects of the
category of locally convex topological vector spaces.
Dually, it is shown (Theorem \ref{T:complementedinj}) that
a complemented subspace of an uncountable direct sum of
Banach spaces is topologically isomorphic to the direct sum
of ({\bf LB})-spaces. This result is used to characterize
(Theorem \ref{T:proj}) projective objects of the above category.} 
\end{abstract}

\maketitle
\markboth{A.~Chigogidze}{Continuous homomorphisms of
Arens-Michael algebras}

\section{Introduction}\label{S:intro}
Arens-Michael algebras are limits of projective systems of Banach algebras
(or, alternatively, closed subalgebras of (uncountable) products
of Banach algebras). Quite often, when dealing with a
particular Arens-Michael algebra, at least one projective system
arises naturally (for instance, as a result of certain construction) and,
in most cases, it does contain the needed information about its limit.
The situation is somewhat different if an Arens-Michael algebra is
given arbitrarily and there is no particular projective system associated
with it in a canonical way.

Below (Definition \ref{D:Frechet}) we introduce the concept
of a projective Fr\'{e}chet system and show (Theorem \ref{T:spectral})
that every continuous homomorphism of Arens-Michael agebras can
be obtained as the limit homomorphism of certain morphism of cofinal
subsystems of the corresponding Fr\'{e}chet systems. This result
applied to the identity homomorphism obviously implies that any
Arens-Michael algebra has essentially unique Fr\'{e}chet
system associated with it. Consequently, any information
about an Arens-Michael algebra is contained in the associated
Fr\'{e}chet system. The remaining problem of restoring this 
information is, of course, still non-trivial, but sometimes can be
successfully handled by using a simple but 
effective method. This method is based on Proposition
\ref{P:3.1.3}\footnote{Applications of Proposition \ref{P:3.1.3} in a
variety of situations can be found in \cite{book}.}. 

Applying such an approach we obtain the following, in a sense
dual, statements (precise statements are recorded in Theorems
\ref{T:complementedpro} and \ref{T:complementedinj} respectively).

\begin{C}
A retract of an uncountable product of Banach spaces
is the product of retracts of countable subproducts.
\end{C}

\begin{D}
A retract of an uncountable coproduct of Banach spaces
is the coproduct of retracts of countable subcoproducts.
\end{D}

Based on these results we present a complete description of
injective and projective objects of the category ${\mathcal LCS}$ of locally
convex topological vector spaces (over the field of complex numbers). Obviously
products (coproducts) of arbitrary collections of injective
(projective) objects of any category are again injective (projective). Also the
class of injectives (projectives) is stable under retractions. Sometimes
in a relatively nice category there naturally exists a class of ``simple"
injectives (projectives) so that any other injective (projective) can
be obtained from simple ones by applying the above mentioned operations,
i.e. by forming products (coproducts) and by passing to retracts.
Consider some examples:
\begin{enumerate}
\item
Injectives in the category ${\mathcal COMP}$ of compact Hausdorff spaces 
 are precisely retracts of products of copies of the
closed unit segment \cite{book}, whereas projectives are retracts of
coproducts of copies of the singleton \cite[Problem 6.3.19(a)]{engelking}.
The cone of the uncountable power of the segment is an example of an
injective which is not the product of simpler injectives.
\item
Injectives in the category ${\mathcal COMPABGR}$ of compact
Hausdorff abelian topological groups  are precisely the tori,
i.e. products of copies of the circle group ${\mathbb T}$ (see \cite{tori}
for a topological characterization of arbitrary tori).
\item
Projectives in the category ${\mathcal R}$-${\mathcal MOD}$ of
left ${\mathcal R}$-modules
(where ${\mathcal R}$ is an associative ring with unit) are 
coproducts of countably generated projectives \cite{kapl}.
\end{enumerate}

In these examples injectives are retracts of products and projectives are
retracts of 
coproducts of simpler objects. On the other hand, as shows the
first example, injectives need not be products.

The following two results (Theorems \ref{T:inj} and \ref{T:proj} respectively)
provide the full solution of the corresponding problems
in the category ${\mathcal LCS}$.

\begin{A}
The following conditions are equivalent for a
locally convex topological vector space $X$:
\begin{itemize}
\item[(1)]
$X$ is an injective object of the category $\mathcal LCS$.
\item[(2)]
$X$ is isomorphic to the product
$\displaystyle \prod\{ F_{t} \colon t \in T\}$,
where each $F_{t}$, $t \in T$,
is a complemented subspace of the product
$\displaystyle \prod\{ \ell_{\infty}(J_{t_{n}}) \colon n \in \omega\}$.
\end{itemize}
\end{A}

%%%%%%%%%%%%%%%%%%%%%%%%%%%%%%%%%%%%%%%%%%%

\section{Preliminaries}\label{S:pre}

\subsection{Projective systems and their morphisms}\label{SS:morphism}

Below we consider projective systems \spep consisting of
topological algebras $X_{\alpha}$, $\alpha \in A$, and 
continuous homomorphisms
$p_{\alpha}^{\beta} \colon X_{\beta} \to X_{\alpha}$,
$\alpha \leq \beta$, $\alpha, \beta \in A$
($A$ is the directed indexing set of ${\mathcal S}_{X}$).
The limit $\varprojlim{\mathcal S}_{X}$ of this system is
defined as the closed subalgebra of the Cartesian product
$\displaystyle \prod\{ X_{\alpha}\colon \alpha \in A\}$
(with coordinatewise defined operations)
consisting of all
{\em threads}
of ${\mathcal S}_{X}$, i.e. 
\[ \varprojlim{\mathcal S}_{X} = \Big\{ x_{\alpha} \in 
\prod\{ X_{\alpha}\colon \alpha \in A\}
\colon
p_{\alpha}^{\beta}(x_{\beta}) = x_{\alpha}\;\;
\text{for any}\;\; \alpha, \beta \in A \;\;\text{with}\;\;
\alpha \leq \beta \Big\} .\]

The $\alpha$-th limit projection
$p_{\alpha} \colon \varprojlim{\mathcal S}_{X} \to X_{\alpha}$,
$\alpha \in A$, of the system ${\mathcal S}_{X}$ is
the restriction
(onto $\varprojlim{\mathcal S}_{X}$) of the $\alpha$-th
natural projection
$\displaystyle \pi_{\alpha} \colon
\prod\{ X_{\alpha}\colon \alpha \in A\} \to X_{\alpha}$.

If $A^{\prime}$ is a
directed subset of the indexing set $A$, then the subsystem
$\{ X_{\alpha}, p_{\alpha}^{\beta}, A^{\prime}\}$ of
${\mathcal S}_{X}$ is denoted ${\mathcal S}_{X}|A^{\prime}$. We
refer the reader to \cite{robert}, \cite{scha} for general
properties of projective systems. 

Suppose we are given two projective systems \spep and
${\mathcal S}_Y = \{ Y_{\gamma}, q_{\gamma}^{\delta}, B \}$
consisting of topological algebras $X_{\alpha}$,
$\alpha \in A$, and $Y_{\gamma}$, $\gamma \in B$.
A {\em morphism of the system ${\mathcal S}_X$ into
the system ${\mathcal S}_Y$} is a family $\{ \varphi ,
\{ f_{\gamma} \colon \gamma \in B \} \}$, consisting of a
nondecreasing function $\varphi \colon B \rightarrow A$ such
that the set $\varphi (B)$ is cofinal in $A$, and of
continuous homomorphisms $f_{\gamma} \colon X_{\varphi (\gamma )}
\rightarrow Y_{\gamma}$ defined for all $\gamma \in B$
such that\\
$$q^{\delta}_{\gamma}f_{\delta} = f_{\gamma}p^{\varphi
(\delta )}_{\varphi (\gamma )},$$
whenever $\gamma , \delta \in B$ and $\gamma \leq \delta$.
In other words, we require (in the above situation) the
commutativity of the following diagram

\[
\begin{CD}
X_{\varphi (\delta )} @>f_{\delta}>> Y_{\delta}\\
@V{p_{\varphi (\gamma )}^{\varphi (\delta )}}VV @VV
{q^{\delta}_{\gamma}}V\\
X_{\varphi (\gamma )} @>f_{\gamma}>> Y_{\gamma}
\end{CD}
\]

\bigskip   

Any morphism\\
\[\{ \varphi , \{ f_{\gamma} \colon \gamma \in B \} \}
\colon {\mathcal S}_X \rightarrow {\mathcal S}_Y\]

\noindent induces a continuous homomorphism, called {\em the
limit homomorphism of the morphism}\index{limit of morphism}\\
\[\varprojlim~\{ \varphi , \{ f_{\gamma} \colon \gamma \in B \}
\} \colon \varprojlim~{\mathcal S}_X \rightarrow
\varprojlim~{\mathcal S}_Y .\]

\noindent To see this, assign to each thread
$x = \{ x_{\alpha} \colon \alpha \in A \}$ of the
system ${\mathcal S}_X$ the point $y = \{ y_{\gamma}
\colon \gamma \in B \}$ of the product $\displaystyle \prod \{
Y_{\gamma} \colon \gamma \in B \}$ by letting\\
\[ y_{\gamma} = f_{\gamma}(x_{\varphi (\gamma )}),
\gamma \in B .\] 

\noindent It is easily seen that the point
$y = \{ y_{\gamma} \colon \gamma \in B \}$ is in fact a
thread of the system ${\mathcal S}_Y$. Therefore,
assigning to $x = \{ x_{\alpha} \colon \alpha \in A \}
\in \varprojlim~{\mathcal S}_X$ the point $y = \{ y_{\gamma}
\colon \gamma \in B \} \in \varprojlim~{\mathcal S}_Y$, we
define a map $\varprojlim~\{ \varphi , \{ f_{\gamma} \colon
\gamma \in B \} \} \colon \varprojlim~{\mathcal S}_X \rightarrow
\varprojlim~{\mathcal S}_Y$. Straightforward verification
shows that this map is a continuous homomorphism. 

Morphisms of projective systems which arise most
frequently in practice are those defined over the same
indexing set. In this case, the map $\varphi \colon
A \rightarrow A$ of the definition of morphism is
taken to be the identity. Below we shall mostly deal
with such situations and use the following notation: 
$\{ f_{\alpha} \colon X_{\alpha} \rightarrow Y_{\alpha};
\alpha \in A \} \colon {\mathcal S}_X
\rightarrow {\mathcal S}_Y$ or sometimes even a shorter form
$\{ f_{\alpha} \} \colon {\mathcal S}_X
\rightarrow {\mathcal S}_Y$.

\begin{pro}\label{P:1.2.13+}
Let \speq be a projective system and $X$ be a topological algebra.
Suppose that for each $\alpha \in A$ a continuous homomorphism $f_{\alpha}
\colon X \rightarrow Y_{\alpha}$ is given in such a way that
$f_{\alpha} = q_{\alpha}^{\beta}f_{\beta}$
whenever $\alpha , \beta \in A$ and $\alpha \leq \beta$.
Then there exists a natural continuous homomorphism $f \colon X
\rightarrow \varprojlim~{\mathcal S}_Y$ (the diagonal product
$\triangle \{ f_{\alpha} \colon \alpha \in A \}$)
satisfying, for each $\alpha \in A$, the condition
$f_{\alpha} = q_{\alpha}f$.
\end{pro}
\begin{proof}
Indeed, we only have to note that $X$,
together with its identity map $id_X$, forms the
projective system ${\mathcal S}$. So the collection
$\{ f_{\alpha} \colon \alpha \in A \}$ is in fact
a morphism ${\mathcal S} \rightarrow {\mathcal S}_Y$.
The rest follows from the definitions given above.
\end{proof}

%%%%%%%%%%%%%%%%%%%%%%%%%%%%%%%%%%%%%%%%

\subsection{Arens-Michael algebras}
\label{S:AMalgebras}
We recall some definitions \cite{helem}.
A {\em polynormed space} $X$ is a
topological linear
space $X$ furnished with a collection $\{ ||\cdot||_{\nu},
\Lambda \}$ of
seminorms generating the topology of $X$.
This simply means
that the collection
\[ \{ x\in X \colon || x-x_{0}||_{\nu} <
\epsilon\}, x_{0}\in X,
\nu \in \Lambda, \epsilon > 0,\]
forms a subbase of the topology of $X$. A
{\em polynormed algebra}
is a polynormed space $X$ wich admits a
separately continuous multiplication bioperator
$m \colon X \times X \to X$. A {\em multinormed algebra}
is a polynormed algebra such that
$|| xy||_{\nu} \leq ||x||_{\nu}\cdot ||y||_{\nu}$ for each
$\nu \in \Lambda$ and any $(x,y) \in X \times X$. Finally,
an {\em Arens-Michael algebra} is a complete
(and Hausdorff) multinormed
algebra.

The following result \cite[Corollary V.2.19]{helem}
provides a characterization of Arens-Michael algebras.
\begin{thm}\label{T:closed}
The following conditions are equivalent for
a multinormed algebra $X$:
\begin{itemize}
\item[(a)]
$X$ is an Arens-Michael algebra.
\item[(b)]
$X$ is limit of a certain projective system of Banach algebras.
\item[(c)]
$X$ is topologically isomorphic to a closed subalgebra of the
Cartesian product of a certain family of Banach algebras.
\end{itemize}
\end{thm}

%%%%%%%%%%%%%%%%%%%%%%%%%%%%%%%%%%%%%%%%%%%%%%%%%
%%%%%%%%%%%%%%%%%%%%%%%%%%%%%%%%%%%%%%%%%%%%%%%%%

\subsection{Set-theoretical facts}\label{SS:set}
For the reader's convenience we begin by presenting necessary
set-theoretic facts. Their complete proofs can be found
in \cite{book}.

Let $A$  be a partially ordered {\em directed set} (i.e.
for every two elements  $\alpha ,\beta \in A$  there exists
an element  $\gamma \in A$  such that  $\gamma \geq \alpha$ 
and  $\gamma \geq \beta$). We say that a subset
$A_1 \subseteq A$ of $A$ {\em majorates} another subset
$A_2 \subseteq A$ of $A$ if for each element $\alpha_2 \in A_2$
there exists an element $\alpha_1 \in A_1$ such that
$\alpha_1 \geq \alpha_2$. A subset which majorates $A$
is called {\em cofinal} in $A$. A subset of  $A$  is said to
be a {\em chain} if every two elements of it are comparable.
The symbol $\sup B$ , where  $B \subseteq A$, denotes the
lower upper bound of $B$ (if such an element exists in $A$).
Let now $\tau$ be an infinite cardinal number. A subset $B$
of $A$  is said to be $\tau$-{\em closed} in $A$ if for each chain
$C \subseteq B$, with ${\mid}C{\mid} \leq \tau$,  we have
$\sup C \in B$, whenever the element $\sup C$ exists in $A$.
Finally, a directed set $A$ is said to be $\tau$-{\em complete}
if for each chain $C$ of elements of $A$ with
${\mid}C{\mid} \leq \tau$, there exists an element
$\sup C$ in $A$. 

The standard example of a $\tau$-complete set can be obtained
as follows. For an arbitrary set $A$ let $\exp A$ denote, as usual,
the collection of all subsets of $A$. There is a natural partial
order on $\exp A$: $A_1 \geq A_2$ if and only if $A_1 \supseteq A_2$.
With this partial order $\exp A$ becomes a directed set.
If we consider only those subsets of the set $A$ which have
cardinality $\leq \tau$, then the corresponding subcollection
of $\exp A$, denoted by $\exp_{\tau}A$, serves as a basic
example of a $\tau$-complete set.

\begin{pro}\label{P:3.1.1}
Let  $\{ A_{t} : t \in T \}$ be a collection of $\tau$-closed and
cofinal subsets of a $\tau$-complete set $A$. If
$\mid T\mid \leq \tau$, then the intersection
$\cap \{ A_{t}: t \in T \}$ is also cofinal
(in particular, non-empty) and $\tau$-closed in $A$ .
\end{pro}

\begin{cor}\label{C:3.1.2}
For each subset $B$, with  $\mid B \mid \leq \tau$, of a
$\tau$-complete set $A$ there exists an element $\gamma \in A$
such that  $\gamma \geq \beta$  for each  $\beta \in B$ .
\end{cor}

\begin{pro}\label{P:3.1.3}
Let  $A$  be a $\tau$-complete set, 
$L \subseteq A^2$, and suppose that the following three
conditions are satisfied:
\begin{description}
\item[Existence] For each $\alpha \in A$ there exists
$\beta \in A$  such that  $(\alpha ,\beta ) \in L$.
\item[Majorantness] If  $(\alpha ,\beta ) \in L$  and
$\gamma \geq \beta$, then  $(\alpha ,\gamma ) \in L$.
\item[$\tau$-closeness] Let $\{ \alpha_{t} : t \in T \}$
be a chain in $A$ with $\mid T\mid \leq \tau$.
If $(\alpha_{t}, \gamma ) \in L$ for some
$\gamma \in A$ and each $t \in T$, then
$(\alpha ,\gamma ) \in L$ where $\alpha =
\sup \{\alpha_{t} \colon t \in T \}$.
\end{description}
   Then the set of all  $L$-{\em reflexive} elements of 
$A$ (an element $\alpha \in A$ is $L$-reflexive if
$(\alpha ,\alpha ) \in L$)  is cofinal and $\tau$-closed in $A$.
\end{pro}

%%%%%%%%%%%%%%%%%%%%%%%%%%%%%%%%%%%%%%%%%%%%%%%%%%%%
%%%%%%%%%%%%%%%%%%%%%%%%%%%%%%%%%%%%%%%%%%%%%%%%%%%%

\section{Continuous homomorphisms of
Arens-Michael algebras}\label{S:homo}
The following statement is needed in the proof of Theorem
\ref{T:spectral}. In the case when all
$X_{\alpha}$'s are Banach algebras its proof can be extracted from
\cite[Proposition 0.1.9]{helem}
(see also \cite[Proof of Proposition V.1.8]{helem}).

\begin{lem}\label{L:factor}
Let \spep be a projective system consisting of Arens-Micheael
algebras $X_{\alpha}$, $\alpha \in A$, $Y$ be a Banach algebra and
$f \colon \varprojlim {\mathcal S} \to Y$ be a continuous
homomorphism. Then there exist an index $\alpha \in A$ and a
continuous
homomorphism $f_{\alpha} \colon X_{\alpha} \to Y$ such
that $f = f_{\alpha}\circ p_{\alpha}$. 
\end{lem}
\begin{proof}

The continuity of $f$ and the definition of the
topology on $\varprojlim{\mathcal S}$ guarantee that
there exists an index $\alpha \in A$ and an open subset
$V_{\alpha} \subseteq X_{\alpha}$ such that
\begin{equation}\label{EQ:1}
f\left( p_{\alpha}^{-1}(V_{\alpha})\right) \subseteq
\left\{  y \in Y \colon |y|\leq 1 \right\} ,
\end{equation}
where $|\cdot |$ denotes the norm of the Banach space $Y$.

Since $X_{\alpha}$ is an Arens-Michael algebra, $X_{\alpha}$ can be
identified with a closed subalgebra of the product
$\displaystyle \prod\{ B_{t} \colon t \in T\}$ of Banach algebras
$B_{t}$, $t \in T$ (Theorem \ref{T:closed}). Let
$||\cdot ||_{t}$ denote the norm of the
Banach space $B_{t}$, $t \in T$. For $S \subseteq T$, let
$\displaystyle \pi_{S} \colon \prod\{ B_{t} \colon
t \in T\} \to \prod\{ B_{t} \colon t \in S\}$ denote the
natural projection
onto the corresponding subproduct. If
$S \subseteq T$ is a finite subset of $T$,
then $||\{ x_{t} \colon t \in S\} ||_{S} =
\max\{ || x_{t}||_{t} \colon t \in S\}$ for each
$\displaystyle \{ x_{t} \colon t \in S\}
\in \prod\{ B_{t} \colon t \in S\}$.

Since $V_{\alpha}$ is open in $X_{\alpha}$, the definition of the
product topology guarantees the existence of a finite subset
$S \subseteq T$ and of a number $\epsilon > 0$ such that
\begin{equation}\label{EQ:2}
\{ x_{\alpha} \in X_{\alpha} \colon
||\pi_{S}(x_{\alpha})||_{S} \leq \epsilon\} \subseteq V_{\alpha}
\end{equation}

Combining (\ref{EQ:1}) and (\ref{EQ:2}), we have

\begin{equation}\label{EQ:3}
f\left( \{ x \in X \colon ||\pi_{S}\left(
p_{\alpha}(x)\right) ||_{S} \leq \epsilon\}\right) \subseteq
\left\{  y \in Y \colon |y|\leq 1 \right\} .
\end{equation}

It then
follows that if $x\in X$
and $||\pi_{S}(p_{\alpha}(x))||_{S}\leq 1$, then
$||\pi_{S}(p_{\alpha}(\epsilon x))||_{S}=\epsilon
||\pi_{S}(p_{\alpha}(x))||_{S}\leq \epsilon$ and consequently

\begin{equation}\label{EQ:4}
\epsilon|f(x)|=|f(\epsilon x)|\leq 1 ,\; 
\text{i.e.}\; |f(x)|\leq\frac{1}{\epsilon} .
\end{equation}

Since 
\[ \left\|
\pi_{S}\left( p_{\alpha}\left( 
\frac{x}{||\pi_{S}(p_{\alpha}(x))||_{S}}\right)\right)
\right\|  =\frac
{1}{||\pi_{S}(p_{\alpha}(x))||_{S}}\cdot
||\pi_{S}(p_{\alpha}(x))||_{S} = 1 ,\]
we must have (by \ref{EQ:4})

\[ \frac{1}{||\pi_{S}(p_{\alpha}(x))||_{S}}|f(x)|=\left|  f\left( 
\frac{x}{||\pi_{S}(p_{\alpha}(x))||_{S}}\right)
\right|  \leq\frac{1}{\epsilon}\]

\noindent and hence 
\begin{equation}\label{EQ:5}
 |f(x)|\leq\frac{1}{\epsilon}||\pi_{S}(p_{\alpha}(x))||_{S}
\;\;\text{for each}\;\; x \in X .
\end{equation}

Let us now show that the map
\[
h(z) = f\left( p_{\alpha}^{-1}\left( \pi_{S}^{-1}(z)
\cap X_{\alpha}\right)\right) \colon
\pi_{S}(p_{\alpha}(X)) \to Y \]
is well defined. Assuming the contrary,
suppose that for some $z \in \pi_{S}(p_{\alpha}(X))$ there exist
two points
$x_{1}, x_{2} \in p_{\alpha}^{-1}\left( \pi_{S}^{-1}(z)
\cap X_{\alpha}\right)$ such that
$f(x_{1}) \neq f(x_{2})$. Consequently,
$|f(x_{1}-x_{2})| \neq 0$.
On the other hand,
$\pi_{S}(p_{\alpha}(x_{1}-x_{2})) = \pi_{S}(p_{\alpha}(x_{1}))-
\pi_{S}(p_{\alpha}(x_{2}))
= z-z = 0$. Then (\ref{EQ:5}) implies that 
\[ 0 \neq |f(x_{1}-x_{2})|\leq\frac{1}{\epsilon
}||\pi_{S}(p_{\alpha}(x_{1}-x_{2}))||_{S} = 0 .
\]
This contradiction shows that the map
$h \colon \pi_{S}(p_{\alpha}(X)) \to Y$ is indeed well defined.
Note that
\begin{equation}\label{EQ:6}
f = h\circ \pi_{S}\circ p_{\alpha} ,
\end{equation}
which implies that the map $h$ is linear. Next consider points 
$z \in \pi_{S}(p_{\alpha}(X))$ and $x\in X$ such that
$\pi_{S}(p_{\alpha}(x)) = z$. By (\ref{EQ:5}),
\[ |h(z)| = |f(x)| \leq\frac{1}{\epsilon
}||\pi_{S}(p_{\alpha}(x))||_{S} = \frac{1}{\epsilon
}||z||_{S}.\]
This shows that $h$ is bounded and, consequently, continuous.
Next let us show that $h \colon \pi_{S}(p_{\alpha}(X)) \to Y$ is
multiplicative. Let
$(x^{\prime}, y^{\prime}) \in \pi_{S}(p_{\alpha}(X))
\times \pi_{S}(p_{\alpha}(X))$ and
consider a point $(x,y) \in X\times X$ such that
$\pi_{S}(p_{\alpha}(x)) = x^{\prime}$ and
$\pi_{S}(p_{\alpha}(y)) = y^{\prime}$. Then, by (\ref{EQ:6}),
\begin{multline*}
h(x^{\prime}\cdot y^{\prime}) =
h\left( \pi_{S}(p_{\alpha}(x)) \cdot \pi_{S}(p_{\alpha}(y))\right) =
h\left( \pi_{S}(p_{\alpha}(x\cdot y))\right) = f(x\cdot y) =\\
 f(x) \cdot f(y) =
h\left( \pi_{S}(p_{\alpha}(x))\right) \cdot
h\left( \pi_{S}(p_{\alpha}(y))\right) =
h(x^{\prime})\cdot h(y^{\prime}) .
\end{multline*}

Since $(Y,|\cdot|)$ is
complete, $h$ admits the
linear continuous
extension 
\[ g \colon \operatorname{cl}_{\prod\{ B_{t} \colon t \in S\}}
\left(\pi_{S}(p_{\alpha}(X))\right) \to Y .\]

Since the multiplication on
$\displaystyle \prod\{ B_{t} \colon t \in T_{f}\}$ is jointly
continuous, we conclude that $g$ is also multplicative.
Finally, define the map $f_{\alpha}$ as the composition
\[ f_{\alpha} = g\circ \left(\pi_{S}|X_{\alpha}\right)
\colon X_{\alpha} \to Y .\]
Obviously $f_{\alpha}$ is a continuous homomorphism
satisfying the required
equality $f_{\alpha}\circ p_{\alpha} = f$.
\end{proof}

Next we introduce the concept of {\em projective
Fr\'{e}chet system}.
\begin{defin}\label{D:Frechet}
Let $\tau \geq \omega$ be a cardinal number. A projective
system \spep consisting of
topological algebras $X_{\alpha}$ and
continuous homomorphisms $p_{\alpha}^{\beta} \colon
X_{\beta} \to X_{\alpha}$,
$\alpha \leq \beta$, $\alpha ,\beta \in A$, is a {\em $\tau$-system if:}
\begin{enumerate}
\item\label{I:1}
$X_{\alpha}$ is a closed subalgebra of the product of at
most $\tau$ Banach algebras, $\alpha \in A$.
\item\label{I:2}
The indexing set $A$ is $\tau$-complete.
\item\label{I:3}
If $\{ \alpha_{\gamma} \colon \gamma \in \tau\}$ is an 
increasing chain
of elements in $A$ with $\alpha = \sup\{ \alpha_{\gamma}
\colon \gamma
\in \tau\}$, then the diagonal product\footnote{See
Proposition \ref{P:1.2.13+}.}
\[ \triangle\{ p_{\alpha_{\gamma}}^{\alpha} \colon \gamma
\in \tau\}
\colon X_{\alpha} \to \varprojlim\{ X_{\alpha_{\gamma}},
p_{\alpha_{\gamma}}^{\alpha_{\gamma+1}},\tau\}\]
is a topological isomorphism.
\item\label{I:4}
$p_{\alpha}(X)$ is dense in $X_{\alpha}$ for each
$\alpha \in A$.
\end{enumerate}
{\em Fr\'{e}chet} systems are defined as projective $\omega$-systems.
\end{defin}

\begin{pro}\label{P:exists}
Every Arens-Micheal algebra $X$
can be represented as the limit of a projective
Fr\'{e}chet system
${\mathcal S}_{X} = \{ X_{A}, 
p_{A}^{B}, \exp_{\omega}T\}$. Conversely,
the limit of a projective Fr\'{e}chet system is an
Arens-Michael algebra.
\end{pro}
\begin{proof}
By Theorem \ref{T:closed}, $X$ can
be identified
with a closed subalgebra of the product
$\displaystyle \prod\{ X_{t} \colon t \in T\}$ of some
collection of Banach algebras.
If
$|T| \leq \omega$, then $X$
itself is a Fr\'{e}chet algebra and therefore our
statement is trivially true. If $|T| > \omega$, then
consider the set $\exp_{\omega}T$ of all
countable subsets of
$T$. Clearly, $\exp_{\omega}T$ is $\omega$-complete set
(see Subsection \ref{SS:set}).
For each $A \in \exp_{\omega}T$,
let $X_{A} = \operatorname{cl}\pi_{A}(X)$
(closure is taken in $\displaystyle \prod\{ X_{t} \colon t
\in A\}$),
where
\[ \pi_{A} \colon \prod\{ X_{t} \colon t \in T \}
\to \prod\{ X_{t} \colon t \in A\}\]
denotes the natural projection onto the
corresponding subproduct.
Also let $p_{A}^{B} =
\pi_{A}^{B}|X_{B}$
where 
\[ \pi_{A}^{B} \colon
\prod\{ X_{t} \colon t \in B\}\to\prod\{ X_{t}
\colon t \in A\}\]
is the natural projection, $A, B \in
\exp_{\omega}T$, $A \leq B$. The straightforward
verification shows
that ${\mathcal S}_{X} = \{ X_{A}, 
p_{A}^{B}, \exp_{\omega}T\}$ is indeed a projective
Fr\'{e}chet system such that
$\varprojlim{\mathcal S}_{X} = X$.

Conversely, let ${\mathcal S}_{X} = \{ X_{\alpha}, 
p_{\alpha}^{\beta}, A\}$ be a projective Fr\'{e}chet system.
Clearly, $\varprojlim{\mathcal S}_{X}$ can be identified with a closed
subalgebra of the product
$\displaystyle \prod\{ X_{\alpha} \colon \alpha \in A\}$
(see Subsection \ref{SS:morphism}). Each
$X_{\alpha}$, $\alpha \in A$, can
obviously be identified with a closed subalgebra of the product
$\displaystyle \prod\{ B_{n}^{\alpha} \colon n \in T_{\alpha}\}$
of a countable
collection of Banach algebras $B_{n}^{\alpha}$. Then
$\varprojlim{\mathcal S}_{X}$, as a closed subalgebra of the
product
$\displaystyle \prod\left\{ \prod\{B_{n}^{\alpha}
\colon n \in T_{\alpha}\}\colon a \in A\right\}$ is, according to
Theorem \ref{T:closed}, an Arens-Michael algebra.
\end{proof}

The following statement is one of our main results.

\begin{thm}\label{T:spectral}
Let $f \colon \varprojlim{\mathcal S}_{X} \to
\varprojlim{\mathcal S}_{Y}$ be a continuous
homomorphism between the limits of two projective Fr\'{e}chet
systems\;\;
\spep\- and\; \speq with the same indexing set $A$.
Then there exist
a cofinal and $\omega$-closed subset $B_{f}$ of $A$
and a morphism
\[ \{ f_{\alpha} \colon X_{\alpha} \to
Y_{\alpha}, B_{f}\} \colon
{\mathcal S}_{X}|B_{f} \to {\mathcal S}_{Y}|B_{f} ,\]
consisting of continuous homomorphisms
$f_{\alpha} \colon X_{\alpha} \to Y_{\alpha}$, $\alpha \in B_{f}$,
such that $f = \varprojlim\{ f_{\alpha}\colon B_{f}\}$.\\

If, in particular, $\varprojlim{\mathcal S}_{X}$ and
$\varprojlim{\mathcal S}_{Y}$
are topologically isomorphic, then \spep and \speq
contain isomorphic
cofinal and $\omega$-closed subsystems.
\end{thm}
\begin{proof}
We perform the spectral search by means of the
following
relation
\begin{multline*}
L = \{ (\alpha,\beta ) \in A^{2}\colon \alpha \leq
\beta\;\text{and there exists a
continuous homomorphism}\\
 f_{\alpha}^{\beta} \colon X_{\beta} \to Y_{\alpha}\;
\text{such that}\; f_{\alpha}^{\beta}p_{\beta} =
q_{\alpha}f \} .
\end{multline*}

Let us verify the conditions of Proposition \ref{P:3.1.3}.

{\bf Existence Condition.} By assumption, $Y_{\alpha}$ is
a Fr\'{e}chet algebra. Therefore $Y_{\alpha}$ can be
identified with a closed subspace of a countable product
$\displaystyle \prod\{ B_{n} \colon n \in \omega\}$ of Banach algebras.
Let
$\displaystyle \pi_{n} \colon \prod\{ B_{n} \colon n \in
\omega\} \to B_{n}$ denote the
$n$-th natural projection. For each $n \in \omega$, by
Lemma \ref{L:factor}, 
there exist an index $\beta_{n} \in A$ and a continuous
homomorphism
$f_{\beta_{n}} \colon X_{\beta_{n}} \to B_{n}$ such that
$\pi_{n}q_{\alpha}f = f_{\beta_{n}}p_{\beta_{n}}$. By
Corollary \ref{C:3.1.2}, there
exists an index $\beta \in A$ such that $\beta \geq
\beta_{n}$ for
each $n$. Without loss of generality we may assume that
$\beta \geq \alpha$. Let $f_{n} =
f_{\beta_{n}}p^{\beta}_{\beta_{n}}$,
$n \in \omega$.
Next consider the diagonal product
\[ f_{\alpha}^{\beta} = \triangle\{ f_{n} \colon n
\in \omega\} \colon X_{\beta} \to \prod\{ B_{n} \colon n
\in \omega\} .\]
Obviously, $f_{\alpha}^{\beta}p_{\beta} = q_{\alpha}f$.
It only remains to show that $f_{\alpha}^{\beta}(X_{\beta})
\subseteq Y_{\alpha}$. First observe that
$f_{\alpha}^{\beta}(p_{\beta}(X)) \subseteq Y_{\alpha}$.
Indeed, let $x \in X$. Then
\begin{multline*}
f_{\alpha}^{\beta}(p_{\beta}(x))
= \{ f_{n}(p_{\beta}(x))\colon n \in \omega \} = \{ f_{\beta_{n}}(p_{\beta_{n}}^{\beta}(p_{\beta}(x)))\colon n
\in \omega \} =\\
 \{ f_{\beta_{n}}(p_{\beta_{n}}(x))\colon n
\in \omega \} = \{ \pi_{n}(q_{\alpha}(f(x))) \colon n
\in \omega \} = q_{\alpha}(f(x)) \in Y_{\alpha}.
\end{multline*}

Finally,
\begin{multline*}
 f_{\alpha}^{\beta}(X_{\beta}) =
f_{\alpha}^{\beta}(\operatorname{cl}
_{X_{\beta}}p_{\beta}(X)) \subseteq
\operatorname{cl}_{\prod\{ B_{n}\colon n \in \omega\}}
f_{\alpha}^{\beta}(p_{\beta}(X)) \subseteq 
\operatorname{cl}_{\prod\{ B_{n}\colon n \in \omega\}}
Y_{\alpha} = Y_{\alpha}.
\end{multline*}

{\bf Majorantness Condition.} The verification of
this condition is trivial. Indeed, it suffiecs to
consider the composition
$f_{\alpha}^{\gamma} = f_{\alpha}^{\beta}
p_{\beta}^{\gamma}$.

{\bf $\omega$-closeness Condition.} Suppose that
for some countable chain $C = \{ {\alpha}_{n}
\colon n \in \omega \}$
in $A$ with $\alpha = \sup C$ and for some $\beta \in A$
with $\beta \geq \alpha$, the maps
$f_{\alpha_{n}}^{\beta} \colon X_{\beta}
\to Y_{\alpha_{n}}$ have already been defined is such a
way that $f_{\alpha_{n}}^{\beta}p_{\beta} =
q_{\alpha_{n}}f$ for each $n \in \omega$ (in other
words, $({\alpha_{n}},\beta ) \in L$ for each
$n \in \omega$). Next consider the composition
\[ f_{\alpha}^{\beta} = i^{-1}\circ \triangle\{
f_{\alpha_{n}}^{\beta} \colon n \in \omega\} \colon X_{\beta}
\to Y_{\alpha} ,\]
where $i \colon Y_{\alpha} \to \varprojlim
\{ Y_{\alpha_{n}}, q_{\alpha_{n}}^{\alpha_{n+1}}, \omega\}$
is the topological isomorphism indicated in condition \ref{I:3}
of Definition \ref{D:Frechet}. Observe that for each $x \in X$
\begin{multline*}
 f_{\alpha}^{\beta}(p_{\beta}(x)) = i^{-1}(\triangle\{
f_{\alpha_{n}}^{\beta} \colon n \in \omega\})(p_{\beta}(x))
=\\
 i^{-1}(\triangle\{ q_{\alpha_{n}}\colon n \in \omega\} )(f(x))
= i^{-1}(\triangle\{ q_{\alpha_{n}}^{\alpha}\colon n
\in \omega\} )
(q_{\alpha}(f(x))) =\\
(\triangle\{ q_{\alpha_{n}}^{\alpha}\colon n \in \omega\} )^{-1}
(\triangle\{ q_{\alpha_{n}}^{\alpha}\colon n \in \omega\} )
(q_{\alpha}(f(x))) = q_{\alpha}(f(x)) .
\end{multline*}
This shows that $(\alpha ,\beta ) \in L$ and finishes the
verification
of the $\omega$-closeness condition.

Now denote by $B_{f}$ the set of all $L$-reflexive elements in $A$.
By Proposition \ref{P:3.1.3}, $B_{f}$ is a cofinal and
$\omega$-closed subset of $A$. One can easily see that the
$L$-reflexivity of an element $\alpha \in A$ is equivalent to
the existence of a continuous homomorphism
$f_{\alpha} = f_{\alpha}^{\alpha} \colon X_{\alpha} \to Y_{\alpha}$
satisfying the equality $f_{\alpha}p_{\alpha} = q_{\alpha}f$.
Consequently, the collection
$\{ f_{\alpha} \colon \alpha \in B_{f} \}$ is a morphism of
the cofinal and $\omega$-closed subspectrum ${\mathcal S}_{X}|B_{f}$
of the spectrum ${\mathcal S}_X$ into the cofinal and
$\omega$-closed subspectrum ${\mathcal S}_{Y}|B_{f}$ of the spectrum
${\mathcal S}_Y$. It only remains to remark that the
original map $f$ is induced by the constructed morphism.
This finishes the proof of the first part of our Theorem.

The second part of this theorem can be obtained from the
first as follows. Let
$f \colon \varprojlim~{\mathcal S}_X \to
\varprojlim~{\mathcal S}_Y$ be a topological isomorphism.
Denote by $f^{-1} \colon 
\varprojlim~{\mathcal S}_Y \to \varprojlim~{\mathcal S}_X$
the inverse of $f$. By the first part proved above,
there exist a cofinal and $\omega$-closed subset $B_{f}$ of
$A$ and a morphism\\
\[ \{ f_{\alpha} \colon X_{\alpha} \to Y_{\alpha} \colon
\alpha \in B_{f} \} \colon {\mathcal S}_{X}|B_{f} \to
{\mathcal S}_{Y}|B_{f}\]
such that $f = \varprojlim \{ f_{\alpha} \colon \alpha
\in B_{f} \}$.
Similarly, there exist a cofinal and $\omega$-closed subset
$B_{f^{-1}}$
of $A$ and a morphism\\
\[ \{ g_{\alpha} \colon Y_{\alpha} \to X_{\alpha} \colon \alpha
\in B_{f^{-1}} \} \colon {\mathcal S}_{Y}|B_{f^{-1}} \to
{\mathcal S}_{X}|B_{f^{-1}}\]
such that $f^{-1} = \varprojlim \{ g_{\alpha} \colon \alpha
\in B_{f^{-1}} \}$.

By Proposition \ref{P:3.1.1}, the set $B = B_{f} \cap B_{f^{-1}}$
is still cofinal and $\omega$-closed in $A$. Therefore, in order
to complete the proof, it suffices to show that for each
$\alpha \in B$ the map $f_{\alpha} \colon X_{\alpha}
\to Y_{\alpha}$ is a topological isomorphism. Indeed, take a point
$x_{\alpha} \in p_{\alpha}\left( \varprojlim~{\mathcal S}_X\right)
\subseteq X_{\alpha}$. Also choose a point $x \in
\varprojlim~{\mathcal S}_X$ such that
$x_{\alpha} = p_{\alpha}(x)$. Then\\
\[ g_{\alpha}f_{\alpha}(x_{\alpha}) = g_{\alpha}
f_{\alpha}p_{\alpha}(x) = g_{\alpha}q_{\alpha}f(x) =
p_{\alpha}f^{-1}f(x) = p_{\alpha}(x) = x_{\alpha}.\]
This proves that
$g_{\alpha}f_{\alpha}|p_{\alpha}\left(
\varprojlim~{\mathcal S}_X\right) =
\operatorname{id}_{p_{\alpha}\left(
\varprojlim~{\mathcal S}_X\right)}$.
Similar considerations show
that $f_{\alpha}g_{\alpha}|q_{\alpha}\left(
\varprojlim~{\mathcal S}_Y\right) =
\operatorname{id}_{q_{\alpha}\left(
\varprojlim~{\mathcal S}_Y\right)}$ for each $\alpha \in B$.
Since $p_{\alpha}\left( \varprojlim~{\mathcal S}_X\right)$
is dense in $X_{\alpha}$
and $q_{\alpha}\left( \varprojlim~{\mathcal S}_Y\right)$
is dense in $Y_{\alpha}$ (condition \ref{I:4} of Definition
\ref{D:Frechet}),
it follows that
$g_{\alpha}f_{\alpha}|X_{\alpha} = \operatorname{id}_{X_{\alpha}}$
and
$f_{\alpha}g_{\alpha}|Y_{\alpha} = \operatorname{id}_{Y_{\alpha}}$.
It is now clear that
$f_{\alpha}$, $\alpha \in B$, is a topological isomorphism
(whose inverse is
$g_{\alpha}$).
\end{proof}

\begin{rem}\label{R:tau}
A similar statement remains true (with the identical proof)
for projective $\tau$-systems
for any cardinal number $\tau > \omega$.
\end{rem}

\begin{rem}\label{R:spectral}
Theorem \ref{T:spectral} is false for countable projective systems.
Indeed, consider the following two projective sequences
\[ {\mathcal S}_{\text{even}} = \{ {\mathbb C}^{2n},
\pi_{2n}^{2(n+1)}, \omega\}\;\;
\text{and}\;\;{\mathcal S}_{\text{odd}} =
\{ {\mathbb C}^{2n+1}, \pi_{2n+1}^{2(n+1)+1}, \omega\} ,\]
\noindent where
\[ \pi_{2n}^{2(n+1)} \colon
{\mathbb C}^{2n} \times {\mathbb C}^{2} \to {\mathbb C}^{2n}\;\; \text{and}\;\;
\pi_{2n+1}^{2(n+1)+1} \colon {\mathbb C}^{2n+1}
\times {\mathbb C}^{2} \to {\mathbb C}^{2n+1}\]
\noindent denote the
natural projections. Clearly the limits
$\varprojlim{\mathcal S}_{\text{even}}$
and $\varprojlim{\mathcal S}_{\text{odd}}$ of these projective systems
are topologically isomorphic (both are topologically isomorphic to the
countable infinite power ${\mathbb C}^{\omega}$ of ${\mathbb C}$),
but ${\mathcal S}_{\text{even}}$
and ${\mathcal S}_{\text{odd}}$ do not contain isomorphic cofinal
subsystems.
\end{rem}

\begin{cor}\label{C:frechet}
Let \spep be a projective Fr\'{e}chet system. If
$\varprojlim{\mathcal S}$
is a Fr\'{e}chet algebra, then there exists an index
$\alpha \in A$ such that
the $\beta$-th limit projection
$p_{\beta} \colon \varprojlim{\mathcal S} \to X_{\beta}$
is a topological isomorphism for each $\beta \geq \alpha$.
\end{cor}
\begin{proof}
Consider a trivial projective Fr\'{e}chet system
${\mathcal S}^{\prime} = \{ X_{\alpha}, q_{\alpha}^{\beta}, A\}$,
where $X_{\alpha} = \varprojlim{\mathcal S}$ and
$q_{\alpha}^{\beta} = \operatorname{id}_{\varprojlim{\mathcal S}}$
for each $\alpha, \beta \in A$. By Theorem \ref{T:spectral}
(applied to the identity homomorphism
$f = \operatorname{id}_{\varprojlim{\mathcal S}}$), there exist
an index $\alpha \in A$ and a continuous homomorphism
$g_{\alpha} \colon X_{\alpha} \to \varprojlim{\mathcal S}$ such that $\operatorname{id}_{\varprojlim{\mathcal S}} = g_{\alpha}\circ p_{\alpha}$.
Clearly, in this situation,
$p_{\alpha}|\varprojlim{\mathcal S} \colon
\varprojlim{\mathcal S} \to X_{\alpha}$
is an embedding with a closed image. But this image
$p_{\alpha}(\varprojlim{\mathcal S})$
is dense in $X_{\alpha}$ (condition (\ref{I:4}) of
Definition \ref{D:Frechet}). Therefore
$p_{\alpha}$ (and, consequently $p_{\beta}$ for
any $\beta \geq \alpha$) is a topological isomorphism. 
\end{proof}

\begin{cor}\label{C:frechetsub}
Suppose that $X$ is a Fr\'{e}chet subalgebra of an uncountable product
$\displaystyle \prod\{ B_{t} \colon t \in T\}$ of Fr\'{e}chet (Banach) algebras,
then there exists a countable subset $T_{X}$ of the indexing set
$T$ such that the restriction $\pi_{T_{X}}|X \colon X \to \pi_{T_{X}}(X)$
of the natural projection
$\displaystyle \pi_{T_{X}} \colon \prod\{ B_{t} \colon t \in T\}
\to \prod\{ B_{t} \colon t \in T_{X}\}$
is a topological isomorphism.
\end{cor}

%%%%%%%%%%%%%%%%%%%%%%%%%%%%%%%%%%%%%%%%%%%%%%%%%%%%%%

\subsection{Arens-Michael ${\ast}$-algebras}\label{S:involutive}
The concept of projective Fr\'{e}chet system can naturally be
adjusted to handle variety of situations. Below it will always
be completely clear in what content this concept is being used.
Let us consider Arens-Michael $\ast$-algebras, i.e. Arens-Michael
algebras with a continuous involution.
It is known that every such an algebra can be identified with a
closed $\ast$-subalgebra of the product of Banach $\ast$-algebras
(see, for instance, \cite[Proposition V.3.41]{helem}).
Therefore one can obtain an alternative description of such
algebras as limits of projective systems consisting of
Banach $\ast$-algebras and continuous $\ast$-homomorphisms
(compare with Theorem \ref{T:closed}). This, as in Proposition
\ref{P:exists}, leads us to the conclusion recorded in the
following statement. 

\begin{pro}\label{P:existsstar}
Every Arens-Micheal $\ast$-algebra $X$
can be represented as the limit of a projective
Fr\'{e}chet system 
${\mathcal S}_{X} = \{ X_{\alpha}, 
p_{\alpha}^{\beta}, A\}$ consisting of Fr\'{e}chet $\ast$-algebras
$X_{\alpha}$, $\alpha \in A$, and continuous $\ast$-homomorphisms
$p_{\alpha}^{\beta} \colon X_{\beta} \to X_{\alpha}$,
$\alpha \leq \beta$, $\alpha ,\beta \in A$. Conversely,
the limit of any such projective Fr\'{e}chet system is an
Arens-Michael $\ast$-algebra.
\end{pro}

The analog of Theorem \ref{T:spectral} is also true.

\begin{pro}\label{P:spectralinvol}
Let $f \colon \varprojlim{\mathcal S}_{X} \to
\varprojlim{\mathcal S}_{Y}$ be a continuous
$\ast$-homomorphism between the limits of two projective Fr\'{e}chet
systems\;\;
\spep\- and\; \speq , consisting of Fr\'{e}chet $\ast$-algebras
and continuous $\ast$-homo\-mor\-p\-hisms and having the same indexing set $A$.
Then there exist
a cofinal and $\omega$-closed subset $B_{f}$ of $A$
and a morphism
\[ \{ f_{\alpha} \colon X_{\alpha} \to
Y_{\alpha}, B_{f}\} \colon
{\mathcal S}_{X}|B_{f} \to {\mathcal S}_{Y}|B_{f} ,\]
consisting of continuous $\ast$-homomorphisms
$f_{\alpha} \colon X_{\alpha} \to Y_{\alpha}$, $\alpha \in B_{f}$,
such that $f = \varprojlim\{ f_{\alpha}; B_{f}\}$.\\

If, in particular, $\varprojlim{\mathcal S}_{X}$ and
$\varprojlim{\mathcal S}_{Y}$
are topologically $\ast$-isomorphic, then \spep and \speq
contain isomorphic
cofinal and $\omega$-closed subsystems.
\end{pro}
\begin{proof}
By Theorem \ref{T:spectral}, there exists a cofinal and
$\omega$-closed subset $B_{f}$ of $A$ and a morphism
\[ \{ f_{\alpha} \colon X_{\alpha} \to
Y_{\alpha}, B_{f}\} \colon
{\mathcal S}_{X}|B_{f} \to {\mathcal S}_{Y}|B_{f} ,\]
consisting of continuous homomorphisms
$f_{\alpha} \colon X_{\alpha} \to Y_{\alpha}$, $\alpha \in B_{f}$,
such that $f = \varprojlim\{ f_{\alpha}; B_{f}\}$.

Let us show that $f_{\alpha}$, $\alpha \in B_{f}$ is actually a
$\ast$-homomorphism. Indeed, let $x_{\alpha} \in p_{\alpha}(X)$ and 
$x \in X$ such that $p_{\alpha}(x) = x_{\alpha}$. Then
\begin{multline*}
 f_{\alpha}(x_{\alpha}^{\ast}) =
f_{\alpha}\left( p_{\alpha}(x)^{\ast}\right) =
f_{\alpha}(p_{\alpha}(x^{\ast})) = q_{\alpha}(f(x^{\ast})) =
q_{\alpha}\left( f(x)^{\ast}\right) =\\ 
q_{\alpha}(f(x))^{\ast} =
f_{\alpha}(p_{\alpha}(x))^{\ast} = f_{\alpha}(x_{\alpha})^{\ast}.
\end{multline*}
\end{proof}

%%%%%%%%%%%%%%%%%%%%%%%%%%%%%%%%%%%%
%%%%%%%%%%%%%%%%%%%%%%%%%%%%%%%%%%%%%%%%%%%%%%%%%%%%%%%%%

\section{Complemented subalgebras of uncountable
products of Fr\'{e}chet algebras}\label{S:complemented}

In this section we show (Theorem \ref{T:complementedpro}) that
complemented subalgebras of (uncountable) products of Fr\'{e}chet
algebras are products of Fr\'{e}chet algebras. We begin
with the following lemma.

\begin{lem}\label{L:retract}
Let $p \colon X \to Y$ be a surjective continuous homomorphism
of topological algebras and
suppose that $X$ is a closed subalgebra of the product
$Y \times B$, where $B$ is a a topological algebra. Assume also that
there exists a continuous homomorphism
$r \colon Y \times B \to X$ satisfying the following conditions:
\begin{itemize}
\item[(i)]
$pr = \pi_{Y}$, where $\pi_{Y} \colon Y \times B \to Y$
denotes the natural projection;
\item[(ii)]
$r(x) = x$ for each $x \in X$.
\end{itemize}
Then there exists a topological isomorphism
$h \colon X \to Y \times \ker p$
such that $\pi_{Y}h = p$.
\end{lem}
\begin{proof}
If $x \in X$, then 
\begin{multline*}
p\left( x - r(p(x),0)\right) = p(x) -p(r(p(x),0))
\stackrel{\text{by} (i)}{=}\\
 p(x) -\pi_{Y}(p(x),0) =
p(x) - p(x) = 0 .
\end{multline*} 
This shows that the formula
\[ h(x) = \left( p(x), x - r(p(x),0)\right) , x \in X,\]
defines a continuous linear map $h \colon X \to Y \times \ker p$.
Moreover, $h$ is a topological isomorphism
between $X$ and $Y \times \ker p$ considered as topological vector
spaces (to see this observe that the continuous and linear map
$g \colon Y \times \ker p \to X$ defined by letting
$g(y, x) = r(y,0)+x$ for each $(y,x) \in Y \times \ker p$, has
the following properties: $g\circ h = \operatorname{id}_{X}$ and
$h\circ g = \operatorname{id}_{Y \times \ker p}$). We
now show that $h$ 
is an isomorophism of the category of topological algebras as
well.

Let $x_{1}, x_{2} \in X$. We need to show that
$h(x_{1})\cdot h(x_{2}) = h(x_{1}\cdot x_{2})$. Since
$X \subseteq Y \times B$ we can write
$x_{i} = (a_{i},b_{i})$, where
$a_{i} \in Y$ and $b_{i} \in B$, $i = 1,2$.
Observe that since $x_{i} \in X$ it follows from (ii)
that $r(x_{i}) = x_{i}$. Consequently, by (i),
$p(x_{i}) = p(r(x_{i})) = \pi_{Y}(x_{i}) =
\pi_{Y}(a_{i},b_{i}) = a_{i}$.
Then 
\[ h(x_{i}) = \left( p(x_{i}), x_{i} - r(p(x_{i}),0)\right) =
\left( a_{i}, (a_{i},b_{i}) -r(a_{i},0)\right) , i = 1,2.\]

\noindent Consequently,
\begin{multline*}
h(x_{1})\cdot h(x_{2}) = \left( a_{1},
(a_{1},b_{1}) -r(a_{1},0)\right) \cdot \left( a_{2},
(a_{2},b_{2}) -r(a_{2},0)\right) = \\
\left( a_{1}\cdot a_{2}, \left[ (a_{1},
b_{1}) -r(a_{1},0)\right] \cdot
\left[ (a_{2},b_{2})-r(a_{2},0)\right]\right) =\\
 \left( a_{1}\cdot a_{2}, \left[ r(a_{1},
b_{1}) -r(a_{1},0)\right] \cdot
\left[ r(a_{2},b_{2})-r(a_{2},0)\right]\right) =\\
\left( a_{1}\cdot a_{2}, r(0,b_{1})\cdot r(0,b_{2})\right)
= \left( a_{1}\cdot a_{2},
r\left[ (0,b_{1})\cdot (0,b_{2})\right]\right) =\\
\left( a_{1}\cdot a_{2}, r(0,b_{1}\cdot b_{2})\right) =
\left( a_{1}\cdot a_{2},r\left[ (a_{1}\cdot a_{2},
b_{1}\cdot b_{2}) - (a_{1}\cdot a_{2},0)\right]\right) = \\
\left( a_{1}\cdot a_{2}, r(a_{1}\cdot a_{2},
b_{1}\cdot b_{2}) - r(a_{1}\cdot a_{2},0)\right) = \\
\left( a_{1}\cdot a_{2}, (a_{1}\cdot a_{2},
b_{1}\cdot b_{2}) - r(a_{1}\cdot a_{2},0)\right) =
h(a_{1}\cdot a_{2}, b_{1}\cdot b_{2}) =\\
 h\left( (a_{1},b_{1})\cdot (a_{2},b_{2})\right) =
h(x_{1}\cdot x_{2}).
\end{multline*}
This shows that $h$ is a homomorphism and, consequently, a
topological isomorphism as required.
\end{proof}

\begin{thm}\label{T:complementedpro}
A complemented subalgebra of the product of uncountable family
of Fr\'{e}chet algebras is topologically isomorphic to the product
of Fr\'{e}chet algebras.
More formally, if $X$ is a complemented subalgebra of the product
$\displaystyle \prod\{ B_{t} \colon t \in T\}$ of
Fr\'{e}chet algebras $B_{t}$, $t \in T$, then $X$ is topologically
isomorphic to the product $\displaystyle \prod\{ F_{j} \colon j \in J\}$,
where $F_{j}$ is a complemented subalgebra of the product
$\displaystyle \prod\{ B_{t} \colon t \in T_{j}\}$ with $|T_{j}| = \omega$
for each $j \in J$.
\end{thm}
\begin{proof}
Let $X$ be a complemented subalgebra of the uncountable product
$\displaystyle B = \prod\{ B_{t} \colon t \in T\}$
of Fr\'{e}chet algebras $B_{t}$, $t \in T$, where $T$ is an
indexing set with $|T| = \tau >\omega$.
There exists a continuous
homomorphism $r \colon B \to X$ such that $r(x) = x$ for each
$x \in X$. A subset $S \subseteq T$ is called $r$-admissible if
$\pi_{S}\left( r(z)\right) = \pi_{S}(z)$ for each point
$z \in \pi_{S}^{-1}\left(\pi_{S}(X)\right)$. 

{\bf Claim 1}. {\em The union of an arbitrary family of
$r$-admissible sets is $r$-admissible.}

Let $\{ S_{j} \colon j \in J\}$ be a collection of
$r$-admissible
sets and $S = \bigcup\{ S_{j} \colon j \in J\}$.
Let $z \in \pi_{S}^{-1}\left( \pi_{S}(X)\right)$.
Clearly $z \in \pi_{S_{j}}^{-1}\left(\pi_{S_{j}}(X)\right)$
for each $j \in J$ and consequently
$\pi_{S_{j}}\left( r(z)\right) = \pi_{S_{j}}(z)$
for each $j \in J$.
Obviously, $\pi_{S}(z) \in \pi_{S}\left( r(z)\right)$
and it therefore suffices to show that the set
$\pi_{S}\left( r(z)\right)$ contains only one point.
Assuming that there is a point
$y \in \pi_{S}\left( r(z)\right)$ such that
$y \neq \pi_{S}(z)$ we conclude (remembering
that $S = \bigcup\{ S_{j} \colon j \in J\}$)
that there must exist an index $j \in J$ such that
$\pi_{S_{j}}^{S}(y) \neq \pi_{S_{j}}^{S}\left(\pi_{S}(z)\right)$.
But this is impossible
\[ \pi_{S_{j}}^{S}(y) \in \pi_{S_{j}}^{S}\left(
\pi_{S}\left( r(z)\right)\right) = \pi_{S_{j}}
\left( r(z)\right) = \pi_{S_{j}}(z) =
\pi_{S_{j}}^{S}\left(\pi_{S}(z)\right) .\]

{\bf Claim 2.} {\em If $S \subseteq T$ is $r$-admissible, then
$\pi_{S}(X)$ is a closed subalgebra of
$\displaystyle B_{S} = \prod\{ B_{t} \colon t \in S\}$.}

Indeed, let $i_{S} \colon B_{S} \to B$ be the
canonical section of $\pi_{S}$
(this means that $i_{S} = \operatorname{id}_{B_{S}}\triangle
\mathbf{0} \colon B_{S} \to B_{S} \times B_{T-S} = B$).
Consider a continuous linear map 
$r_{S} = \pi_{S}\circ r\circ i_{S} \colon B_{S} \to \pi_{S}(X)$.
Obviously, $i_{S}(y) \in \pi_{S}^{-1}\left(\pi_{S}(X)\right)$
for any point $y \in \pi_{S}(X)$. Since $S$ is $r$-admissible
the latter implies that
\[ y = \pi_{S}\left(i_{S}(y)\right) =
\pi_{S}\left( r\left(i_{S}(y)\right)\right) = r_{S}(y) .\]
This shows that $\pi_{S}(X)$ is closed in $B_{S}$.

{\bf Claim 3}. {\em Let $S$ and $R$ be $r$-admissible subsets
of $T$ and $S \subseteq R \subseteq T$. Then the map
$\pi_{S}^{R} \colon \pi_{R}(X) \to \pi_{S}(X)$
is topologically isomorphic to the natural projection
$\pi \colon \pi_{S}(X) \times \operatorname{ker}
\left(\pi_{S}^{R}|\pi_{R}(X)\right) \to \pi_{S}(X)$.}

Obviously $\pi_{R}(X) \subseteq \pi_{S}(X) \times B_{R-S}
\subseteq B_{R} = B_{S} \times B_{R-S}$. Consider the map
$i_{R} = \operatorname{id}_{B_{R}}\triangle
\mathbf{0} \colon B_{R} \to B_{R} \times B_{T-R} = B$. Also
let $r_{R} = \pi_{R}\circ r\circ i_{R} \colon B_{R} \to \pi_{R}(X)$.

Observe that
$\pi_{S}^{R}\circ r_{R}|\left( \pi_{S}(X) \times B_{R-S}\right) =
\pi_{S}^{R}|\left( \pi_{S}(X) \times B_{R-S}\right)$. Indeed, if
$x \in  \pi_{S}(X) \times B_{R-S}$, then $i_{R}(x)
\in \pi_{S}^{-1}\left(\pi_{S}(X)\right)$. Since $S$ is $r$-admissible,
we have $\pi_{S}\left(r\left(i_{R}(x)\right)\right) =
\pi_{S}\left(i_{R}(x)\right)$. Consequently,
\begin{multline*}
\pi_{S}^{R}\left( r_{R}(x)\right) =
\pi_{S}^{R}\left(\pi_{R}\left( r\left(i_{R}(x)\right)\right)\right) = \pi_{S}\left(r\left(i_{R}(x)\right)\right) =
\pi_{S}\left(i_{R}(x)\right) =\\
\pi_{S}^{R}\left(\pi_{R}\left( i_{R}(x)\right)\right) =
\pi_{S}^{R}(x) .
\end{multline*}

Next observe that $r_{R}(x) = x$ for any point $x \in \pi_{R}(X)$.
Indeed, since $R$ is $r$-admissible and since
$i_{R}(x) \in \pi_{R}^{-1}\left(\pi_{R}(X)\right)$ we have
\[ r_{R}(x) = \pi_{R}\left(r\left(i_{R}(x)\right)\right) =
\pi_{R}\left(i_{R}(x)\right) = x .
\]

Application of Lemma \ref{L:retract} (with $X = \pi_{R}(X)$,
$Y = \pi_{S}(X)$, $B = B_{R-S}$, $p = \pi_{S}^{R}|\pi_{R}(X)$ and
$r = r_{R}$) finishes the proof of Claim 3.

{\bf Claim 4.} {\em Every countable subset of $T$ is
contained in a countable $r$-admissible subset of $T$.}

Let $A$ be a countable subset of $T$. Our goal is to find a countable
$r$-admissible subset $C$ such that $A \subseteq C$. By Theorem
\ref{T:spectral}, there exist a countable subset $C$ of $T$
and a continuous homomorphism $r_{C} \colon B_{C} \to \operatorname{cl}_{B_{C}}\left(\pi_{C}(X)\right)$ such
that $A \subseteq C$ and $\pi_{C}\circ r = r_{C}\circ \pi_{C}$.
Consider a point $y \in \pi_{C}(X)$. Also pick a point $x \in X$
such that $\pi_{C}(x) = y$. Then
\[ y  = \pi_{C}(x) = \pi_{C}\left( r(x)\right) =
r_{C}\left(\pi_{C}(x)\right) = r_{C}(y) .\]
This shows that $r_{C}|\pi_{C}(X) = \operatorname{id}_{\pi_{C}(X)}$.
It also follows that $\pi_{C}(X)$ is closed in $B_{C}$.

In order to show that $C$ is $r$-admissible let us consider a point
$z \in \pi_{C}^{-1}\left(\pi_{C}(X)\right)$. By the observation
made above,
$r_{C}\left(\pi_{C}(z)\right) = \pi_{C}(z)$. Finally
\[ \pi_{C}(z) = r_{C}\left(\pi_{C}(z)\right) =
\pi_{C}\left( r(x)\right) \]
which implies that $C$ is $r$-admissible.

Since $|T| = \tau$, we can write
$T = \{ t_{\alpha} \colon \alpha < \tau\}$.
Since the collection of countable $r$-admissible
subsets of $T$ is cofinal in
$\exp_{\omega}T$ (see Claim 4), each element
$t_{\alpha} \in T$ is contained in a countable $r$-admissible
subset $A_{\alpha} \subseteq T$. According to Claim 1, the set
$T_{\alpha} = \bigcup\{ A_{\beta} \colon
\beta \leq \alpha\}$ is $r$-admissible for each $\alpha < \tau$.
Consider the projective system
\[ {\mathcal S} = \{ X_{\alpha}, p_{\alpha}^{\alpha +1}, \tau\} ,\]
where 
\[ X_{\alpha} = \pi_{T_{\alpha}}(X)\;\;\text{and}\;\;
p_{\alpha}^{\alpha +1} = \pi_{T_{\alpha}}^{T_{\alpha +1}}|
\pi_{T_{\alpha +1}}(X)\;\; \text{for each}\;\; \alpha < \tau .\]
Since $T = \bigcup\{ T_{\alpha} \colon \alpha < \tau\}$, it follows
that $X = \projlim\mathcal S$. Obvious transfinite induction based on
Claim 3 shows that 
\[ X = \projlim{\mathcal S} = X_{0} \times\prod\{
\operatorname{ker}\left( p_{\alpha}^{\alpha +1}\right) \colon
\alpha < \tau \} .\]
Since, by the construction, $A_{\alpha}$ is a countable
$r$-admissible subset of $T$, it follows from Claim 2 that
$X_{0}$ and $\ker\left( p_{\alpha}^{\alpha +1}\right)$, $\alpha <\tau$,
are Fr\'{e}chet algebras. This finishes the proof of
Theorem \ref{T:complementedpro}.
\end{proof}
%%%%%%%%%%%%%%%%%%%%%%%%%%%
%%%%%%%%%%%%%%%%%%%%%%%%%%%

\section{Injective objects of the category ${\mathcal LCS}$}
\label{S:injectives}

In this section we investigate injective objects of the category
$\mathcal LCS$ of locally convex Hausdorff topological
vector spaces and their continuous linear maps.

Recall that an object $X$ of the category
${\mathcal LSC}$ is {\em injective} if any continuous linear map
$f \colon A \to X$, defined on a linear subspace of a space $B$, admits
a continuous linear extension $g \colon B \to X$ (i.e. $g|A = f$).
Here is the corresponding diagram

\begin{picture}(300,120)
\put(101,100){$B$}
\put(100,10){$A$}
\put(106,23){\vector(0,1){73}}
\put(99,20){$\cup$}
\put(80,55){$\operatorname{incl}$}
\put(118,12){\vector(1,0){105}}
\put(227,10){$X$}
\put(110,97){\vector(3,-2){115}}
\put(160,18){$f$}
\put(160,68){$g$}
\end{picture}

We start with the metrizable case. The following statement is probably
known (consult with \cite[Theorem C.3.4]{pietsch}, \cite[pp.6--9]{cigler}).

\begin{pro}\label{P:binj}
The following conditions are equivalent for any Banach space
$X$:
\begin{itemize}
\item[(1)]
$X$ is an injective object of the category $\mathcal LCS$.
\item[(2)]
$X$ is an injective object of the category $\mathcal BAN$.
\item[(3)]
$X$ is isomorphic to a complemented subspace of
$\ell_{\infty}(J)$ for some set $J$.
\end{itemize}
\end{pro}
\begin{proof}
The implication $(1) \Longrightarrow (2)$ is trivial.

$(2) \Longrightarrow (3)$. The Banach space $X$ can be identified
with a closed linear subspace of the space $\ell_{\infty}(J)$ for
some set $J$. By condition (2), there exists a linear continuous map
$r \colon \ell_{\infty}(J) \to X$ such that $r(x) = x$ for each
$x \in X$. This obviously implies that $X$ is a complemented
subspace of $\ell_{\infty}(J)$.

$(3) \Longrightarrow (1)$. First let us show that the Banach
space $\ell_{\infty}(J)$ (for any set $J$) is an injective
object of the category ${\mathcal BAN}$. Let
\[ f \colon X \to \ell_{\infty}(J) = \left\{ \{ x_{j} \colon j \in J\}
\in \prod\{ {\mathbb C}_{j} \colon j \in J\} \colon
\sup\{ ||x_{j}|| \colon j \in J\} < \infty \right\} ,\]
where ${\mathbb C}_{j}$, $j \in J$, stands for a copy of ${\mathbb C}$, be a
continuous linear map defined on a closed linear subspace $X$ of a
Banach space $Y$. Since $f$ is bounded, for each $j \in J$ we have
\[ ||f(x)_{j}|| \leq \sup\{ ||f(x)_{j}|| \colon j \in J\} =
||f(x)|| \leq ||f||\cdot ||x|| , \; x \in X .\] 
This shows that $||\pi_{j}\circ f|| \leq ||f||$, where
$\pi_{j} \colon \ell_{\infty}(J) \to {\mathbb C}_{j}$ denotes the
canonical projection onto the $j$-th coordinate. By the
Hahn-Banach Theorem,
the linear map $\pi_{j}\circ f \colon X \to {\mathbb C}_{j}$, $j \in J$,
admits a continuous linear
extension $g_{j} \colon Y \to {\mathbb C}_{j}$ such that
$||g_{j}|| = ||f\circ \pi_{j}|| \leq ||f||$. Consequently,
\[\sup\{ ||g_{j}(y)|| \colon j \in J\} \leq
\sup\{ ||g_{j}||\cdot ||y|| \colon j \in J\} \leq
||f||\cdot ||y|| < \infty\; ,\; y \in Y .\]
This shows that the map
$g \colon Y \to \ell_{\infty}(J)$, given by letting
\[ g(y) = \{ g_{j}(y) \colon j \in J\} , y \in Y,\]
is well defined. Obviously $g$ is continuous, linear and
extends the map $f$. Therefore $\ell_{\infty}(J)$ is indeed an injective
object of the category ${\mathcal BAN}$.

Next consider a complete locally convex
topological vector space $Y$,
its closed linear subspace $Z$ and a continuous linear map
$f \colon Z \to X$, where $X$ is a complemented subspace
of the Banach space $\ell_{\infty}(J)$ for some set $J$.
We need to show that $f$ admits a continuous
linear extension $g \colon Y \to X$. Identify $Y$ with a
closed linear subspace of the product
$\displaystyle \prod\{ B_{t} \colon t \in T\}$ of
Banach spaces $B_{t}$, $t \in T$. Clearly $Z$ is a closed
linear subspace of $\displaystyle \prod\{ B_{t} \colon t \in T\}$.
By Lemma \ref{L:factor},
there exist a finite subset $T_{f}$ and
a continuous linear map $f^{\prime} \colon
\operatorname{cl}_{B_{T_{f}}}\left( \pi_{T_{f}}(Z)\right) \to X$
such that $f = f^{\prime}\circ\pi_{T_{f}}|Z$. Clearly
$\operatorname{cl}_{B_{T_{f}}}\left( \pi_{T_{f}}(Z)\right)$
is a closed linear subspace of the Banach space
$\operatorname{cl}_{B_{T_{f}}}\left( \pi_{T_{f}}(Y)\right)$.
The first part of the proof of this implication,
coupled with condition $(3)$, implies that
$X$ is an injective object of the category $\mathcal BAN$. Therefore
the
map $f^{\prime}$ admits a continuous linear extension
$g^{\prime} \colon \operatorname{cl}_{B_{T_{f}}}
\left( \pi_{T_{f}}(Y)\right) \to X$. It only remains to note that
the map $g = g^{\prime}\circ\pi_{T_{f}}|Y$ is a continuous
linear extension of $f$.
\end{proof}

\begin{thm}\label{T:inj}
The following conditions are equivalent for a
locally convex topological vector space $X$:
\begin{itemize}
\item[(1)]
$X$ is an injective object of the category $\mathcal LCS$.
\item[(2)]
$X$ is isomorphic to the product
$\displaystyle \prod\{ F_{t} \colon t \in T\}$,
where each $F_{t}$, $t \in T$,
is a complemented subspace of the product
$\displaystyle \prod\{ \ell_{\infty}(J_{t_{n}}) \colon n \in \omega\}$.
\end{itemize}
\end{thm}
\begin{proof}
$(2) \Longrightarrow (1)$. By Proposition \ref{P:binj},
$\ell_{\infty}(J)$ is an injective
object of the category $\mathcal LCS$ for any set $J$. Obviously
the product of an
arbitrary collection of injective objects of the category
$\mathcal LCS$ is also an injective
object of this category. Consequently, the Fr\'{e}chet space
$F_{t}$, $t \in T$, as a complemented subspace of
$\displaystyle \prod\{ \ell_{\infty}(J_{t_{n}})\colon n \in \omega\}$
is injective. Finally, the space $X$, as a product of injectives,
is an injective object of the category $\mathcal LCS$.

$(1) \Longrightarrow (2)$. The space $X$ can be identified
with a closed linear subspace of the product
$\displaystyle \prod\{ B_{t} \colon t \in T\}$ of
Banach spaces $B_{t}$, $T \in T$. Each of the spaces $B_{t}$ can in turn
be identified with a closed linear subspace of the space
$\ell_{\infty}(J_{t})$ for some set $J_{t}$, $t \in T$.
Condition (1) implies in this situation that $X$
is a complemented subspace of the product 
$\displaystyle\prod\{ \ell_{\infty}(J_{t}) \colon t \in T\}$.
The required conclusion now follows from Theorem
\ref{T:complementedpro}.
\end{proof}

 %%%%%%%%%%%%%%%%%%%%%%%%%%%%%%%%%%%%%%%%%%%%%%%%%%%%%%%
%%%%%%%%%%%%%%%%%%%%%%%%%%%%%%%%%%%%%%%%%%%%%%%%%%%%%

\section{Projective objects of the category ${\mathcal LCS}$}
\label{S:projectives}
\subsection{Complemented subspaces of uncountable sums of Banach spaces}\label{SS:complementedinj}
In this section we consider coproducts in the
category of complete locally convex topological vector spaces.
These are called locally convex direct sums (\cite[p. 55]{scha}
or simply direct sums \cite[p. 89]{robert}. We are mainly interested
in direct sums $\displaystyle \bigoplus\{ B_{t} \colon t \in T\}$ of
(uncountable) collections of Banach spaces $B_{t}$, $t \in T$.
Let us recall corresponding definitions. Let
$\{ B_{t} \colon t \in T\}$
be an arbitrary collection of (Banach) spaces. For each
$t \in T$ the space
$B_{t}$ is identified (via an obvious isomorphism) with
the following subspace
\[ B_{t} = \left\{ \{ x_{t} \colon t \in T\} \in \prod\{
B_{t} \colon t \in T\} \colon x_{t^{\prime}} = 0\;
\text{whenever}\; t^{\prime}\neq t\right\}\]
of the product $\displaystyle\prod\{ B_{t} \colon t \in T\}$ of (Banach)
spaces $B_{t}$, $t \in T$. As a set the direct sum
$\displaystyle \bigoplus\{ B_{t} \colon t \in T\}$ coincides with
the subset of the product $\displaystyle\prod\{ B_{t} \colon t \in T\}$
consisting of points with finitely many non-zero coordinates
(in other words, $\displaystyle \bigoplus\{ B_{t} \colon t \in T\}$ is the
vector subspace of $\displaystyle\prod\{ B_{t} \colon t \in T\}$
spanned by
$\displaystyle\bigcup\{ B_{t} \colon t \in T\}$). The topology on
$\displaystyle \bigoplus\{ B_{t} \colon t \in T\}$ is the finest
locally convex topology for which each of the natural embeddings
$\displaystyle B_{t} \hookrightarrow \bigoplus\{ B_{t} \colon t \in T\}$
is continuous. Categorical description of this construction is
also useful. Let $B$ be a locally convex topological vector space
containing each of the space $B_{t}$, $t \in T$, as a subspace.
Suppose that for any choice of continuous linear maps
$f_{t} \colon B_{t} \to C$, where $C$ is a locally convex
topological vector space, there exists an unique continuous linear
map $f \colon B \to C$ such that $f|B_{t} = f_{t}$
for each $t \in T$. Then $B$ is canonically isomorphic
to the direct sum $\displaystyle \bigoplus\{ B_{t} \colon t \in T\}$.
A straightforward verification shows that if $S \subseteq T$,
then $\displaystyle \bigoplus\{ B_{t} \colon t \in S\}$ can be canonically
identified with the subspace of $\displaystyle \bigoplus\{ B_{t} \colon t \in T\}$
consisting of points $\displaystyle \left\{ \{ x_{t} \colon t \in T\} \in
\bigoplus\{ B_{t} \colon t \in T\} \colon x_{t} = 0\;
\text{for each}\; t \in T-S\right\}$. Note also that if
$S \subseteq R \subseteq T$, then the map
\[ \pi_{S}^{R} \colon \bigoplus\{ B_{t} \colon t \in R\} \to
\bigoplus\{ B_{t} \colon t \in S\} ,\]
defined by letting\label{PAGE:?}
\[
\pi_{S}^{R}\left(\{ x_{t} \colon t \in R\}\right) = 
\begin{cases}
x_{t} ,\;\text{if}\; t \in S\\
0\; , \;\;\text{if}\; t \in R-S ,
\end{cases}
\]
\noindent is continuous and linear.

\begin{pro}\label{P:injspectral}
Let $\displaystyle f \colon \bigoplus\{ B_{t} \colon t \in T\}
\to \bigoplus\{ B_{t} \colon t \in T\}$ be a
continuous linear map of the direct sum of
Banach spaces
$B_{t}$, $t \in T$, into itself. Suppose that $S \subseteq T$ and
$\displaystyle f\left(\bigoplus\{ B_{t} \colon t \in S\}\right)
\subseteq \bigoplus\{ B_{t} \colon t \in S\}$.
If $|T| >\omega$, then the collection
\[ {\mathcal K}_S = \left\{ A \in \exp_{\omega}(T-S)
\colon f\left(\bigoplus\{ B_{t} \colon t \in S \cup A\}\right)
\subseteq \bigoplus\{ B_{t} \colon t \in S \cup A\}\right\} \]
\noindent is cofinal and $\omega$-closed in $\exp_{\omega}(T-S)$.
\end{pro}
\begin{proof}
Consider the following relation $L_S$ on the
set $(exp_{\omega}(T-S))^2$:\\
\begin{multline*}
L_S = \left\{ (A,C) \in (exp_{\omega}(T-S))^2
\colon A \subseteq C \; \text{and}\right.\\
\left. f\left(\bigoplus\{ B_{t} \colon t \in S \cup A\}\right) \subseteq
\bigoplus\{ B_{t} \colon t \in S \cup C\}\right\} 
\end{multline*}
\noindent Next we perform the
spectral search with respect to $L_S$.

{\bf Existence Condition.} We have to show
that for each $A \in exp_{\omega}(T-S)$ there
exists $C \in exp_{\omega}(T-S)$ such that $(A,C) \in L_S$.

We begin with the following observation.

{\bf Claim.} {\em For each $j \in T$ there exists a
finite subset $C_{j} \subseteq T$ such that
$\displaystyle f\left( B_{j}\right) \subseteq
\bigoplus\{ B_{t} \colon t \in C_{j}\}$.}

The unit ball $X = \{ x \in B_{j} \colon
||x||_{j} \leq 1\}$ (here $||\cdot ||_{j}$ denotes a
norm of the Banach space $B_{j}$) being bounded in $B_{j}$
is, by \cite[Theorem 6.3]{scha}, bounded in
$\displaystyle \bigoplus\{ B_{t} \colon t \in T\}$. Continuity of $f$
guarantees that $f(X)$ is also bounded in
$\displaystyle \bigoplus\{ B_{t} \colon t \in T\}$.
Applying \cite[Theorem 6.3]{scha} once again, we
conclude that there exists
a finite subset $C_{j} \subseteq T$ such that
$\displaystyle f(X) \subseteq \bigoplus\{ B_{t} \colon t \in C_{j}\}$.
Finally the linearity of $f$ implies that
$\displaystyle f\left( B_{j}\right) \subseteq
\bigoplus\{ B_{t} \colon t \in C_{j}\}$.

Let now $A \in \exp_{\omega}(T-S)$. For each $j \in A$,
according to Claim,
there exists a finite subset $C_{j} \subseteq T$ such that
$\displaystyle f\left( B_{j}\right) \subseteq
\bigoplus\{ B_{t} \colon t \in C_{j}\}$.
Let $\displaystyle \widehat{C} = \cup\{ C_{j} \colon j \in A\}$.
Observe that $|\widehat{C}| \leq \omega$.
Clearly $\displaystyle f\left( B_{j}\right) \subseteq \bigoplus\{ B_{t}
\colon t \in \widehat{C}\}$ for each $j \in A$. The linearity of $f$
guarantees in this situation that
\[ f\left( \bigoplus\{ B_{t} \colon t \in A \}\right)
\subseteq \bigoplus\{ B_{t} \colon t \in \widehat{C}\} .\]
Finally let $C = (\widehat{C}-S) \cup A$. Since, by our assumption,
$\displaystyle f\left(\bigoplus\{ B_{t} \colon t \in S\}\right)
\subseteq \bigoplus\{ B_{t} \colon t \in S\}$, it follows that
$\displaystyle f\left(\bigoplus\{ B_{t} \colon t \in S \cup A\}\right)
\subseteq \bigoplus\{ B_{t} \colon t \in S \cup C\}$.
Therefore $(A,C) \in L_{S}$.

{\bf Majorantness Condition.} Let $(A,C) \in L_S$,
$D \in \exp_{\omega}(T-S)$ and $C \subseteq D$.
Condition $(A,C) \in L_S$ implies that
$\displaystyle f\left(\bigoplus\{ B_{t} \colon t \in S \cup A\}\right)
\subseteq \bigoplus\{ B_{t} \colon t \in S \cup C\}$. The inclusion
$C \subseteq D$ implies that
$\displaystyle \bigoplus\{ B_{t} \colon t \in C\}
\subseteq \bigoplus\{ B_{t} \colon t \in D\}$. Consequently
\[ f\left(\bigoplus\{ B_{t} \colon t \in S \cup A\}\right)
\subseteq \bigoplus\{ B_{t} \colon t \in S \cup C\} \subseteq
\bigoplus\{ B_{t} \colon t \in S \cup D\} ,\]
which means that $(A,D) \in L_{S}$.

{\bf $\omega$-closeness Condition.} Let
$(A_i ,C) \in L_S$, $i \in \omega$, where
$\{ A_i \colon i \in \omega \}$ is a countable
chain in $\exp_{\omega}(T-S)$. We need to
show that $(A,C) \in L_S$, where
$A = \cup \{ A_i \colon i \in \omega \}$. 
Condition $(A_{i},C) \in L_{S}$ implies that
$\displaystyle f\left(\bigoplus\{ B_{t} \colon t
\in S \cup A_{i}\}\right) \subseteq \bigoplus\{ B_{t}
\colon t \in S \cup C\}$, $i \in \omega$. Also observe that
\[ \bigoplus\left\{ \bigoplus\{ B_{t} \colon t \in A_{i}\} \colon i
\in \omega\right\} = \bigoplus\{ B_{t} \colon t \in A\} .\]
Consequently, by the linearity of $f$,
\begin{multline*}
f\left( \bigoplus\{ B_{t} \colon t \in S \cup A\}\right) =
f\left(  \bigoplus\left\{ \bigoplus\{ B_{t} \colon t
\in S \cup A_{i}\} \colon i
\in \omega\right\}\right) \subseteq\\
 \bigoplus\{ B_{t} \colon t \in S \cup C\} .
\end{multline*}

By Proposition \ref{P:3.1.3}, the set of
$L_{S}$-reflexive elements
of $\exp_{\omega}(T-S)$
is cofinal and $\omega$-closed
in $\exp_{\omega}(T-S)$. It only remains to note that
an element $A \in \exp_{\omega}(T-S)$ is
$L_{S}$-reflexive if and only if $A \in {\mathcal K}_S$.
\end{proof}

\begin{rem}\label{R:i}
Proposition \ref{P:injspectral} is valid only for uncountable
sums (compare with Remark \ref{R:spectral}).
In order to see this consider an uncountable
sum $\displaystyle \bigoplus\{ B_{t} \colon t \in T\}$ of Banach spaces and
suppose that
$\displaystyle f \colon \bigoplus\{ B_{t} \colon t \in T\}
\to \bigoplus\{ B_{t} \colon t \in T\}$ is a topological isomorphism.
Proposition
\ref{P:injspectral}, applied to the map $f$, guarantees the existence of
a cofinal and $\omega$-closed subset
${\mathcal K}_{\emptyset}(f)$ of $\exp_{\omega}T$
satisfying the following property: $\displaystyle
f\left(\bigoplus\{ B_{t} \colon t \in A\}\right) \subseteq
\left(\bigoplus\{ B_{t} \colon t \in A\}\right)$
for each $A \in {\mathcal K}_{\emptyset}(f)$. The same Proposition
\ref{P:injspectral}, applied to the map $f^{-1}$ (recall that $f$ is a
topological isomorphism), guarantees the existence of a cofinal
and $\omega$-closed subset
${\mathcal K}_{\emptyset}(f^{-1})$ of $\exp_{\omega}T$
satisfying the following property:
$\displaystyle f^{-1}\left(\bigoplus\{ B_{t}
\colon t \in A\}\right) \subseteq
\left(\bigoplus\{ B_{t} \colon t \in A\}\right)$
for each $A \in {\mathcal K}_{\emptyset}(f^{-1})$. By Proposition
\ref{P:3.1.1}, the intersection
${\mathcal K} = {\mathcal K}_{\emptyset}(f) \cap
{\mathcal K}_{\emptyset}(f^{-1})$
is still cofinal and $\omega$-closed in $\exp_{\omega}T$.
Obviously, for
each $A \in {\mathcal K}$ the restriction 
\[ f\left|\bigoplus\{ B_{t} \colon t \in A\}\right.
\colon \bigoplus\{ B_{t} \colon t \in A\} \to
\bigoplus\{ B_{t} \colon t \in A\}\]
is a topological isomorphism.

Again, for countable sums such a phenomenom is impossible.
In order to see
this consider the inductive sequences
\[ {\mathcal S}_{\text{even}} =
\{ {\mathbb C}^{2n}, i_{2n}^{2(n+1)}, \omega\}\;\;
\text{and}\;\; {\mathcal S}_{\text{odd}} =
\{ {\mathbb C}^{2n+1}, i_{2n+1}^{2(n+1)+1}, \omega\} ,\] 
\noindent where
\[  i_{2n}^{2(n+1)} \colon
{\mathbb C}^{2n} \hookrightarrow {\mathbb C}^{2n} \bigoplus
{\mathbb C}^{2}\;\text{and}\;\;
 i_{2n+1}^{2(n+1)+1} \colon {\mathbb C}^{2n+1}
\hookrightarrow {\mathbb C}^{2n+1} \bigoplus {\mathbb C}^{2}\]
\noindent denote the
natural embeddings. Clearly the limit spaces
$\varinjlim{\mathcal S}_{\text{even}}$
and $\varinjlim{\mathcal S}_{\text{odd}}$ of these inductive systems
are topologically isomorphic, but ${\mathcal S}_{\text{even}}$
and ${\mathcal S}_{\text{odd}}$ do not contain isomorphic cofinal
subsystems.
\end{rem}

\begin{thm}\label{T:complementedinj}
A complemented subspace of the sum of uncountable family
of Banach spaces is isomorphic to the sum of $({\mathbf L\mathbf B})$-spaces.
More formally, if $X$ is a complemented subspace of the sum
$\displaystyle \bigoplus\{ B_{t} \colon t \in T\}$ of
Banach spaces $B_{t}$, $t \in T$, then $X$ is
isomorphic to the sum $\displaystyle \bigoplus\{ F_{j} \colon j \in J\}$
where $F_{j}$ is a complemented subspace of the countable sum
$\displaystyle \bigoplus\{ B_{t} \colon t \in T_{j}\}$ where $|T_{j}| = \omega$
for each $j \in J$.
\end{thm}
\begin{proof}
Let $X$ be a complemented subspace of an uncountable sum
$\displaystyle \bigoplus\{ B_{t} \colon t \in T\}$ of
Banach spaces $B_{t}$, $t \in T$. Clearly there exists a continuous
linear map $\displaystyle r \colon \bigoplus\{ B_{t} \colon t \in T\} \to X$
such that $r(x) = x$ for each $x \in X$. A subset $S \subseteq T$
is called $r$-admissible if
$\displaystyle r\left(\bigoplus\{ B_{t} \colon t \in S\}\right) \subseteq
\bigoplus\{ B_{t} \colon t \in S\}$.

For a subset $S \subseteq T$,
let $\displaystyle X_{S} = r\left(\bigoplus
\{ B_{t} \colon t \in S\}\right)$.

We state properties of $r$-admissible subsets needed below.

{\bf Claim 1}.
{\em If $S \subseteq T$ is an $r$-admissible, then
$\displaystyle X_{S} = X \bigcap \left(\bigoplus\{ B_{t}
\colon t \in S\}\right)$.}

Indeed, if $y \in X_{S}$, then there exists a point
$\displaystyle x \in \bigoplus\{ B_{t} \colon t \in S\}$ such that
$r(x) = y$. Since $S$ is $r$-admissible, it follows that
\[\displaystyle y = r(x) \in r\left(\bigoplus
\{ B_{t} \colon t \in S\}\right) \subseteq
\bigoplus\{ B_{t} \colon t \in S\} .\]
\noindent Clearly, $y \in X$. This shows that
$\displaystyle X_{S} \subseteq X \bigcap \left(\bigoplus\{ B_{t}
\colon t \in S\}\right)$.

Conversely. If
$\displaystyle y \in X \bigcap \left(\bigoplus\{ B_{t}
\colon t \in S\}\right)$, then $y \in X$ and hence,
by the property of $r$,
$y = r(y)$. Since $\displaystyle y \in \bigoplus\{ B_{t}
\colon t \in S\}$, it follows that
$\displaystyle y = r(y) \in r\left(\bigoplus\{ B_{t}
\colon t \in S\}\right) = X_{S}$.

{\bf Claim 2}.
{\em The union of an arbitrary collection of $r$-admissible
subsets of $T$ is $r$-admissible.}

This is matter of a straightforward
verification of the definition of the $r$-ad\-missi\-bi\-lity.

{\bf Claim 3.}
{\em Every countable subset of $T$ is contained in a countable
$r$-admissible subset
of $T$.}

This follows from Proposition \ref{P:injspectral}
applied to the map $r$.

{\bf Claim 4.}
{\em If $S\subseteq T$ is an $r$-admissible subset of $T$, then
$r_{S}(x) = x$ for each point $x \in X_{S}$, where 
$\displaystyle r_{S} = r\left|\left( \bigoplus\{
B_{t} \colon t \in S\}\right.\right) \colon \bigoplus\{
B_{t} \colon t \in S\} \to X_{S}$.}

This follows from the corresponding property of the map $r$.

{\bf Claim 5.}
{\em If $S$ and $R$ are $r$-admissible subset of $T$ and $S \subseteq R$,
then $X_{R}$ is isomorphic to the sum
$\displaystyle X_{S}\bigoplus \ker \left( f_{S}^{R}\right)$,
where the map $f_{S}^{R} \colon X_{R} \to X_{S}$ is defined by
letting\footnote{Recall that the definition of the map $\pi_{S}^{R}$ is given
right before Proposition  \ref{P:injspectral} on page \pageref{PAGE:?}.}
$f_{S}^{R}(x) = r_{S}\left( \pi_{S}^{R}(x)\right)$ for each $x \in X_{R}$.}

Observe that $f_{S}^{R}(x) = x$ for each $x \in X_{S}$. 

{\bf Claim 6.} 
{\em If $S$ and $R$ are $r$-admissible subsets of $T$ and $S \subseteq R$,
then
there exists a continuous linear map
$i_{S}^{R} \colon \ker \left( f_{S}^{R}\right) \to
\ker\left( \pi_{S}^{R}\right)$  such that
$r_{R} \circ i_{S}^{R} = \operatorname{id}_{\ker\left( f_{S}^{R}\right)}$}.

Let $x \in \ker \left( f_{S}^{R}\right)$. Then
$r_{S}\left(\pi_{S}^{R}(x)\right) = 0$ and consequently,
\[\pi_{S}^{R}\left( x-\pi_{S}^{R}(x)\right) =
\pi_{S}^{R}(x)-\pi_{S}^{R}\left(\pi_{S}^{R}(x)\right) =
\pi_{S}^{R}(x) - \pi_{S}^{R}(x) = 0 .\]
This shows that by letting $i_{S}^{R}(x) = x-\pi_{S}^{R}(x)$ for each
$x \in \ker\left( f_{S}^{R}\right)$, we indeed define
a map $i_{S}^{R} \colon \ker \left( f_{S}^{R}\right) \to
\ker\left( \pi_{S}^{R}\right)$. Finally observe that
\begin{multline*}
 r_{S}^{R}\left( i_{S}^{R}(x)\right) =
r_{R}\left( x-\pi_{S}^{R}(x)\right) =
r_{R}(x) - r_{R}\left(\pi_{S}^{R}(x)\right) =
r_{R}(x) - r_{S}\left(\pi_{S}^{R}(x)\right)=\\
 r_{R}(x) = x .
\end{multline*}
In other words, $r_{S}^{R}\circ i_{S}^{R} =
\operatorname{id}_{\ker\left( f_{S}^{R}\right)}$
as required.

Let $|T| = \tau$. Then we can write
$T = \{ t_{\alpha} \colon \alpha < \tau\}$.
Since the collection of countable $r$-admissible
subsets of $T$ is cofinal in
$\exp_{\omega}T$ (see Claim 3), each element
$t_{\alpha} \in T$ is contained in a countable $r$-admissible
subset $A_{\alpha} \subseteq T$. According to Claim 2, the set
$T_{\alpha} = \bigcup\{ A_{\beta} \colon
\beta \leq \alpha\}$ is $r$-admissible for each $\alpha < \tau$.
Consider the inductive system
\[ {\mathcal S} = \{ X_{\alpha}, i_{\alpha}^{\alpha +1}, \tau\} ,\]
where 
\[ X_{\alpha} = X_{T_{\alpha}} = X\bigcap r\left(\bigoplus\{ B_{t}
\colon t \in T_{\alpha}\}\right)\;\;\text{(see Claim 1)}\]
and
\[ i_{\alpha}^{\alpha +1} \colon X_{\alpha} \to X_{\alpha +1}\]
denotes the natural inclusion for each $\alpha < \tau$.
For a limit ordinal number $\beta < \tau$ the sum
$\displaystyle \bigoplus\{ B_{t} \colon t \in T_{\beta}\}$ is
topologically isomorphic to the limit space of the direct system
$\displaystyle \left\{ \bigoplus\{ B_{t} \colon t \in T_{\alpha}\},
j_{T_{\alpha}}^{T_{\alpha +1}}, \alpha < \beta \right\}$, where
\[\displaystyle j_{T_{\alpha}}^{T_{\alpha +1}} \colon
\bigoplus\{ B_{t} \colon t \in T_{\alpha}\} \to
\bigoplus\{ B_{t} \colon t \in T_{\alpha +1}\}\]
\noindent is the natural inclusion. This observation, coupled with Claim 4,
implies that
$X_{\beta} = \injlim\{ X_{\alpha}, i_{\alpha}^{\alpha +1}, \alpha <\beta\}$
for each limit ordinal number $\beta < \tau$.

In particular $X$ is topologically isomorphic to the limit space
of the inductive system
$\{ X_{\alpha}, i_{\alpha}^{\alpha +1}, \alpha <\tau\}$.

For each $\alpha < \tau$, according to Claim 5, the inclusion
$i_{\alpha}^{\alpha +1} \colon X_{\alpha} \to X_{\alpha +}$ is
topologically isomorphic to the inclusion
$\displaystyle X_{\alpha} \hookrightarrow X_{\alpha} \bigoplus
\ker\left( f_{T_{\alpha}}^{T{\alpha +1}}\right)$. In this situation the
straightforward transfinite induction shows that $X$ is topologically
isomorphic to the direct sum 
\[\displaystyle X_{0}\bigoplus \left(\bigoplus\left\{
\ker\left( f_{T_{\alpha}}^{T_{\alpha +1}}\right) \colon
\alpha < \tau \right\}\right) .\]

Since, by construction, the set $T_{0}$ is countable, Claim 4
guarantees that $X_{0}$ is an $({\bf LB})$-space. Similarly, since the
set $T_{\alpha +1}-T_{\alpha} = A_{\alpha +1}$ is countable,
Claim 6 guarantees
that the space
$\displaystyle \ker\left( f_{T_{\alpha}}^{T_{\alpha +1}}\right)$
is an $({\bf LB})$-space for each $\alpha < \tau$. This completes the proof
of Theorem \ref{T:complementedinj}.
\end{proof}

%%%%%%%%%%%%%%%%%%%%%%%%%%%%%%%%%%%%%%
%%%%%%%%%%%%%%%%%%%%%%%%%%%%%%%%%%%%%%

%%%%%%%%%%%%%%%%%%%%%%%%%%%%%

\end{document}